\documentclass[ reqno,final]{article}
\usepackage{amssymb,amsmath, amsthm,latexsym,graphicx}
\usepackage[color]{showkeys}
\usepackage{hyperref}
\definecolor{refkey}{gray}{.75}
\usepackage{listings}
\usepackage{booktabs}

\hypersetup{linktocpage}
\hypersetup{
    colorlinks,
    citecolor=black,
    filecolor=black,
    linkcolor=black,
    urlcolor=black
}

\newtheorem{theorem}{Theorem}[section]

\theoremstyle{remark}

\def\clap#1{\hbox to 0pt{\hss#1\hss}}

\def\p{\partial}
\def\tilde{\widetilde}
\def\bs{\boldsymbol}

\def\ud{\,\mathrm{d}}

\def\ub{\bs{u}}

\title{Goal-oriented A Posteriori Error Estimation for Finite Volume Methods}
\author{Qingshan Chen\thanks{Corresponding author: Department of
    Scientific Computing, Florida State University, Tallahassee, FL
    32306. Email: {\tt qchen3@fsu.edu}, Url:
    \url{http://people.sc.fsu.edu/\~qchen3}}\;  
    \and Max Gunzburger} 

\date{ }

\begin{document}

\maketitle

\begin{abstract}
A general framework for goal-oriented a posteriori error
  estimation for 
  finite volume methods is presented. The framework does not rely on
  recasting finite volume methods as special cases of finite element
  methods, but instead directly determines error estimators from the
  discretized finite volume equations.   Thus, the framework can be
  applied to 
  arbitrary finite volume methods. It also provides the proper
  functional settings to address well-posedness issues for the
  primal and adjoint problems. Numerical results are presented to
  illustrate the validity and effectiveness of the a posteriori error
  estimates and their applicability to adaptive mesh refinement.
\end{abstract}
%


\section{Introduction}\label{sec:intro}

Finite volume methods have become increasingly popular due to
their intrinsic conservative properties and their capability in dealing
with complex domains; see, e.g.~\cite{FiPy:2009,RTKS10}. Hence, a
posteriori error estimates for finite volume methods are important as
they aid in  
error control and improve the overall accuracy of numerical
simulations. A posteriori error estimates also play a key role in the
implementation of adaptive mesh refinement methods. 
In fact, the main motivation of the current work is to find simple
and robust a posteriori error estimators to guide adaptive mesh
refinements for finite volume methods in regional climate modeling
\cite{CGR, RJGJDS}.  

The literature on a posteriori error analysis for finite volume
methods is slim 
compared to that for finite element methods. A 
large part of the available works aim to derive a posteriori error
bounds for approximate solutions in certain global energy norms;
see, e.g.~\cite{ABC03,ABMV03,AgOu00, CaLT05, Nic05, Nic06, Voh08a,
  Voh08}.  

We are particularly interested in another type of a 
posteriori error estimates that are goal oriented, i.e.,~estimates of
errors in certain quantities of interest. A goal-oriented error
estimate is potentially very useful in assisting in error
control. However, the literature on goal-oriented a posteriori error
estimates 
for finite volume methods is even scarcer, probably due to the fact that
finite volume methods do not naturally fit into variational
frameworks. Insightful efforts have been made to address this challenge
by exploiting the equivalence between certain finite volume methods
and numerical schemes in variational forms such as the finite
element methods or discontinuous Galerkin methods. For example, in
\cite{BarLar02}, a goal-oriented a posteriori
error 
estimate is presented for a special type of finite volume method that is
equivalent to a Petrov-Galerkin variant of the discontinuous
Galerkin method. In \cite{EPPTW09}, an a
posteriori error analysis is presented for cell-centered finite volume
methods for the convection-diffusion problem by utilizing the 
equivalence between the finite volume methods and the lowest-order
Raviart-Thomas mixed finite element method with a special quadrature.
However, the applicability of this approach is limited for two
reasons. First, on many grids on which finite volume schemes are
constructed, there are no quadrature rules known, and thus
implementing finite element or discontinuous Galerkin schemes on these
grids is impractical. An example of such grids is hexagonal
Voronoi grid (\cite{DFGun99}). The other reason is that finite volume methods
for real-world problems are often sophisticated in themselves, and it is
often not clear, to say the least, how to establish a connection to 
schemes in variational forms. In this regard we again refer to  \cite{RTKS10}.

In this work, we aim to derive a general functional analytic framework
for a posteriori error estimation for arbitrary finite volume methods. 
The idea is to derive a posteriori error estimators at
the partial differential equation level in an appropriate functional setting. 
This approach does not require the differential equations or the
numerical schemes to be recast in variational form nor do they rely on
connecting an finite volume method with a finite element method. 
 The approximate
solutions produced by finite volume methods are simply taken as inputs
to the a posteriori error estimator. Because the a posteriori error estimation is
independent of the exact form of the finite volume method, it can be
applied to arbitrary finite volume methods. 

In Galerkin finite element or discontinuous Galerkin methods, the
difference $u-u_h$ between the exact solution $u$ and the approximate
solution $u_h$ is orthogonal to the test function space $V_h$. For
this reason, a posteriori error estimation for these methods
requires that the adjoint solution $\phi$ be sought in a space
$V_{h'}$ larger than $V_h$. This restriction does not apply in our
approach due to the very fact that finite volume schemes do not
naturally fit into variational forms, though for the sake of accuracy,
it may be advantageous to seek the adjoint solution with a higher
order scheme. This point will be made clear in the next section. 

It has been pointed out that the well-posedness issue for the adjoint
equation remains challenging and open in many cases; see, e.g.,
\cite{BarLar02}.  
A byproduct of our approach is that, because the
adjoint problem is naturally posed in an appropriate functional
setting, its well posedness can be dealt with by the
abundant analytical tools of standard partial differential equation
theories; again, see \cite{BarLar02}.

The rest of the paper is organized as follows. At the beginning of the 
next section, we present our approach for  a posteriori error
estimation for finite volume methods in a general functional
analytical framework. It is 
followed by a numerical example demonstrating the error
estimates. In Section \ref{sec:appl-adapt-mesh}, we present an
application of the a posteriori error estimation to adaptive mesh
refinement. We conclude with some remarks in Section
\ref{sec:conclusion}. 

\section{A posteriori error estimates for finite volume
  methods}\label{sec:a-post-error-estim} 
\subsection{Abstract framework}\label{sec:framework}

Let $H$ denote a Hilbert space endowed with the norm $||\cdot||$ and 
inner product $(\cdot ,\cdot)$. Let $L$ denote an unbounded operator in
$H$ with domain $D(L)$ dense in $H$. 
The primal problem  we  deal with is succinctly
formulated in this functional setting as:\\
\indent{\em for each $f\in H$, find $\ub\in D(L)$ such that}
\begin{equation}
  \label{eq:2}
  L\ub = f.
\end{equation}
We consider time-dependent
problems so that the operator $L$ usually takes the form 
\begin{displaymath}
  L = \dfrac{\p}{\p t} + A,
\end{displaymath}
where $A$ represents a linear differential spatial operator that is
usually also unbounded in a respective function space. 
 In this work, we assume that the primal problem \eqref{eq:2}
is well-posed, i.e.,~it possesses a unique solution.

The quantity of interest is given as a possibly nonlinear functional $Q(\ub)$ 
 of $\ub$. 
 Let $\ub^\#$ denote an approximate solution of the primal problem \eqref{eq:2}.
Then, the error in the quantity of interest can be written as 
   \begin{equation}
\label{eq:3}
       Q(\ub) - Q(\ub^\#) = (\ub-\ub^\#, \phi)
   \end{equation}
   for some kernel function $\phi\in H$. We note that the kernel function
   $\phi$ may depend on $\ub$; indeed this is the case when $Q$ is a
   nonlinear function of $u$; see Section \ref{sec:appl-adapt-mesh}.

The foregoing assumption that the domain $D(L)$ is dense in $H$ is
crucial as it allows us to rigorously define the adjoint operator
$L^*$ of $L$ and its domain $D(L^*)$. Indeed, according to
\cite{Ru91}, 
a function $\tilde\ub\in H$ belongs to $D(L^*)$ if and only if 
the mapping $D(L): \ub \longrightarrow (L\ub,\,\tilde\ub)$ 
is a linear bounded functional on $D(L)$ for the norm of $H$. Because
$D(L)$ is dense in $H$, the Hahn-Banach theorem guarantees that the bounded linear
functional can be extended to the whole of $H$. Then, by
the Riesz representation theorem, there exists a unique element,
denoted as $L^*\tilde\ub$, of $H$ such that  
\begin{equation}
  \label{eq:4}
  (L\ub,\,\tilde\ub) = (\ub,\,L^*\tilde\ub)\qquad \forall\,\ub \in
  D(L). 
\end{equation}
Thus, $L^*$ is a linear, possibly unbounded, operator from $D(L^*)$ to
$H$. 
The adjoint problem can be formally stated as:\\
\indent{\em for a given function $\phi\in H$, find $\tilde\ub\in
  D(L^*)$ such that }
\begin{equation}
  \label{eq:5}
  L^*\tilde\ub = \phi.
\end{equation}

We will demonstrate, through examples, how to compute  adjoint
operators. Generally speaking, for time-dependent problems, the
operator $L^*$ usually has the form of   
\begin{displaymath}
  L^* = -\dfrac{\p}{\p t} - B. 
\end{displaymath}
Primal time-dependent problems are usually
{\em initial}-boundary value problems. Consequently, most adjoint
problems are {\em final}-boundary value problem because the
values of the unknown function are
imposed at the final time $t= T$. Applying the change of variable
$
  \tau = T -t
$
transforms the final-value problem into an initial-value
problem. In this way, many analytic techniques can be employed to establish the
well-posedness of the adjoint problem. We will demonstrate this point
with specific examples. 

Now that all the necessary functional settings have been introduced,
we shall derive an a posteriori error estimate for the forward problem
\eqref{eq:2}, regardless of the finite volume methods actually used to
solve the primal and adjoint problems.
By \eqref{eq:3}--\eqref{eq:5},
we infer that
\begin{align*}
  Q(\ub) - Q(\ub^\#) =& (\ub-\ub^\#,\,\phi)\\
=& (\ub-\ub^\#,\,L^*\tilde\ub)\\
=& (L(\ub-\ub^\#),\,\tilde\ub)\\
=& (f-L\ub^\#,\,\tilde\ub).
\end{align*}
Therefore, we have
\begin{equation}
  \label{eq:6}
   Q(\ub) - Q(\ub^\#) =  (f-L\ub^\#,\,\tilde\ub).
\end{equation} 
Note that the true solution $\ub$ is not required for computing
the error and that, if the adjoint solution $\tilde\ub$ is exact,
then the error estimate is actually exact as well. These are the
key advantages of the a posteriori estimate \eqref{eq:6}. The cost incurred is, of  
course, the need to solve the adjoint problem \eqref{eq:5}. 

So far we have not discussed the exact form of the
finite volume methods. The approximate solutions produced by these 
schemes are taken as input to \eqref{eq:6}. Hence, the
formula, in principle, applies to arbitrary finite volume methods. 
We should also note that the approximate solution to the primal
problem is usually given in a discrete form and thus the term
$L\ub^\#$ in \eqref{eq:6} is not well defined at this point. We will
explore a walk-around to this issue when we deal with specific
examples in Section \ref{sec:numer-demonstr-1d} and
\ref{sec:appl-adapt-mesh}.

\subsection{Numerical demonstration: one-dimensional scalar equation
}\label{sec:numer-demonstr-1d} 
The primary goal of this section is, by the means of a simple example,
to demonstrate the implementation of the abstract framework laid out
in the previous section. We will also demonstrate the validity of the
a posteriori error estimates by comparing the estimated errors with
the true errors. Finally, we will explore the impact of the numerical errors
in the adjoint solution on a posteriori error estimates, and what
measures can be taken to ensure accuracy. 

We consider the one-dimensional linear transport equation
\begin{align}
   & u_t + au_x = f, & &0<x<1, t>0,\label{e1.1}\\
   & u(0,t) = g(t),  & &t>0,\label{e1.2}\\
   & u(x,0) = u_0(x), & &0<x<1.\label{e1.3}
\end{align}
We assume that the coefficient $a$ is positive and constant, in which
case the initial and boundary value problem
\eqref{e1.1}--\eqref{e1.3} is well-posed, provided that the problem
is cast in an appropriate functional setting. The well-posedness of
this problem will not be discussed here because it is a classical
example in partial differential equation theory; see, among many texts,
\cite{Evans10}. Nevertheless, we shall specify the proper functional
settings for the problem so that its adjoint problem can be defined
and an a posteriori error estimate for its solution can be derived.

We let $H =  L^2\left( (0,1)\times(0,T) \right)$ and let $L$ be the linear
operator associated with \eqref{e1.1}--~\eqref{e1.3}. 
We define the domain $D(L)$ as 
$$
   \mathcal{D}(L) = \left\{u\in H \,\,|\,\, u_t + au_x\in H,\, u(0,t) = 0,\,u(x,0) = 0  \right\}.
$$
It can shown that $D(L)$ is dense in $H$.
For each $u\in \mathcal{D}(L)$, the operator $L$ is defined as
$$
   L u = u_t + au_x.
$$

We now define the adjoint operator $L^*$ of $L$ and its domain $\mathcal{D}(L^*)$.
We recall that the domain $\mathcal{D}(L^*)$ is defined as a space of functions
$\tilde u$ for which
$(Lu, \tilde u)$ is a continuous functional for $u\in\mathcal{D}(L)$
with respect to the norm of $H$, that 
is, 
\begin{displaymath}
   \tilde u\in \mathcal{D}(L^*) \iff u\rightarrow (Lu,\tilde u) \textrm{ is continuous in the $H$ norm.}
\end{displaymath}
Based on this definition, we determine that $\mathcal{D}(L^*)$ is given by
$$
   \mathcal{D}(L^*) = \left\{\tilde u\in H \,\,|\,\, \tilde u_t + a\tilde
     u_x\in H,\, \tilde u(1,t) = 0,\,\tilde u(x,T) = 0  \right\}
$$
and, for each $\tilde u\in \mathcal{D}(L^*)$, we have 
$$
   L^* \tilde u = -\tilde u_t - a\tilde u_x.
$$
For $u\in\mathcal{D}(L)$ and $\tilde u\in\mathcal{D}(L^*)$, the following 
important relationship
holds:
\begin{equation}
   (Lu,\,\tilde u) = (u, \, L^*\tilde u).
   \label{e1.15}
\end{equation}
Within this functional setting, the adjoint problem can be
stated as follows:\\
\indent{\em  for any } $\phi \in H$, {\it find} $\tilde u \in D(L^*)$ {\it such that} 
\begin{equation}
  \label{eq:7}
  L^* \tilde u = \phi.
\end{equation}

We now consider the well-posedness issue for
the adjoint problem \eqref{eq:7}.
The equations of this problem, in differential
form, read
\begin{align}
   -\tilde u_t - a\tilde u_x &= \phi\qquad\textrm{ for }0<x<1,\,\, 0<t<T,\label{e1.21a}\\
   \tilde u(x,T) &= 0\qquad\textrm{ for } 0<x<1,\label{e1.21b}\\
   \tilde u(1,t) &= 0\qquad\textrm{ for } 0<t<T.\label{e1.21c}
\end{align}
The conditions \eqref{e1.21b} and \eqref{e1.21c} are due to the
requirement that $v$  
should be sought in the adjoint domain $\mathcal{D}(L^*)$. 
The system \eqref{e1.21a}--\eqref{e1.21c} needs to be solved
backward in time because a {\it final} condition
is given at $t=T$. To overcome this awkwardness, we make 
the change of variable
$
   t = T-\tau.
$
Then, the system \eqref{e1.21a}--\eqref{e1.21c} becomes 
\begin{align}
   &\tilde u_\tau - a\tilde u_x = \phi,\label{e1.25a}\\
   &\tilde u(x,0) = 0,\label{e1.25b}\\
   &\tilde u(1,\tau) = 0,\label{e1.25c}
\end{align}
which can be solved forwards in time.
The well-posedness of this system can be studied in 
the same way as that for the primal system \eqref{e1.1}--\eqref{e1.3}, e.g.,~by semigroup theory. 
In the language of that theory, the operator
$L^*$ is the infinitesimal generator a semigroup of 
contractions $S(\tau)$.
The well-posedness result for the system above
is stated in the following theorem.

\begin{theorem}
   If $\phi\in L^1(0,T;L^2(0,1))$, then there exists
   an unique solution $\tilde u\in C([0,T];L^2(0,1))$ of 
   \eqref{e1.25a}--\eqref{e1.25c} of the form
   $$
      \tilde u(\tau) = \int_0^\tau S(\tau - s)\phi(s)\ud s.
   $$
\end{theorem}
\noindent We do not prove the theorem here and instead
refer to \cite{Pazy83}.

The simple
example considered in this section demonstrates a key advantage of our approach:
we derive and use an adjoint system whose 
well-posedness can be studied in the same 
manner as that for the primal system,
for which a multitude of analytical tools are available. 

We now concern ourselves with linear quantities of
interest that can be expressed as 
\begin{equation}
  \label{eq:8}
  Q(u) = \int_\Omega u\phi,
\end{equation}
where $\phi$ is a kernel function independent of the unknown $u$.
Identifying $\phi$ in \eqref{eq:8} with $\phi$ in \eqref{e1.25a}, the
a posteriori error estimate of $Q(u) - Q(u^\#)$ is then given by
\eqref{eq:6}. 
However, there is an implementation issue associated with \eqref{eq:6}. In
practice, we only have the approximate solution  $u^\#$ in its discrete
form, e.g., as approximations of the function values at
discrete grid points. For
\eqref{eq:6} to make sense, we could either  define a discrete
analogue ${L}^\#$ of ${L}$ or define and apply a
mapping that maps $u^\#$ from its discrete representation to a
continuous representation so that the differential operator
${L}$ can be applied. In this work, we take a third approach. We 
find a way to 
transpose $L$ back onto $\tilde u$ and then  use 
\eqref{eq:7} to replace $L^*\tilde u$ by $\phi$. This approach
avoids approximating $L u^\#$ with discrete values which usually
introduces another source of error. However, we note
that, in general, $ u^\#$ does not belong to $\mathcal{D}(L)$ so that
we cannot
invoke \eqref{e1.15} for $(L u^\#, v)$. Thus, we proceed by integration by 
parts:
\begin{align*}
   (L u^\#, v) &= \int_I( u^\#_t + a u^\#_x)v\\
   &= -\int_{t=0}u^\#v - a\int_{x=0} u^\#v + (u^\#,\phi).
\end{align*}
Therefore,
\begin{equation}
   Q(u) - Q(u^\#) = (f,v)  - (u^\#,\phi) +
   \int_{t=0}u_0v + a\int_{x=0}gv .
   \label{e1.17}
\end{equation}

\subsubsection{Numerical results}

The forward problem \eqref{e1.1}--\eqref{e1.3} is solved with the
first-order explicit upwind method. The exact form of the scheme is
not essential to this work, but we shall briefly state the method for
the sake of reference. (Other finite volume methods used in the sequel,
however, will not be explicitly described.) We choose $0=x_0 < x_1 < \cdots < x_M =
1$ and let $K_i = [x_{i-1},\,x_i]$ denote the $i^{\textrm{th}}$ control
volume/cell. Obviously, $\{K_i\}_{i=1}^M$ forms a partition of the
interval 
$[0,1]$. We also set $\Delta t = T/N$ and let $t_n = n\Delta t$ denote the
discrete time steps. We let $U^n_i$ denote the average of the unknown
$u$ over the cell $K_i$ at time $t=t^n$. Then, the first-order explicit
upwind method for \eqref{e1.1} can be written as 
$$
  \dfrac{U^{n+1}_i - U^n_i}{\Delta t} + a\dfrac{U^n_i -
    U^n_{i-1}}{\delta_i} = F^n_i \quad\textrm{ for } 1\leq i\leq M.
$$
In the above, $\delta_i = |K_i|$ and $U_0$ is the average of the
unknown $u$ on an artificial 
cell $K_0 = [x_{-1},\,x_0]$ with $x_{-1} = -x_1$. The boundary
condition \eqref{e1.2} can be imposed as 
$$
  \dfrac{1}{2}\left(U^{n+1}_0 + U^{n+1}_1\right) =
  \dfrac{1}{2}g(t^n). 
$$

To demonstrate the validity of the a posteriori error estimate
\eqref{e1.17} for the problem \eqref{e1.1}--\eqref{e1.3},
we consider an example with a sine wave solution:
\begin{align}
   &u_t + u_x = 0,\qquad 0<x<1,\,0<t<0.5,\label{e1a.1}\\
   &u(x,0) = \sin(2\pi x)\label{e1a.2},\\
   &u(0,t) = -\sin(2a\pi t).\label{e1a.3}
\end{align}
This problem has an analytic solution 
\begin{equation}
  \label{eq:9}
  u(x,t) = \sin(2\pi(x-at)),
\end{equation}
which allows us to study the performance of the a posteriori error
estimate \eqref{e1.17} by comparing the estimated errors to the true
errors. 

\begin{figure}[h!]
   \begin{center}
    \includegraphics[width=4.5in]{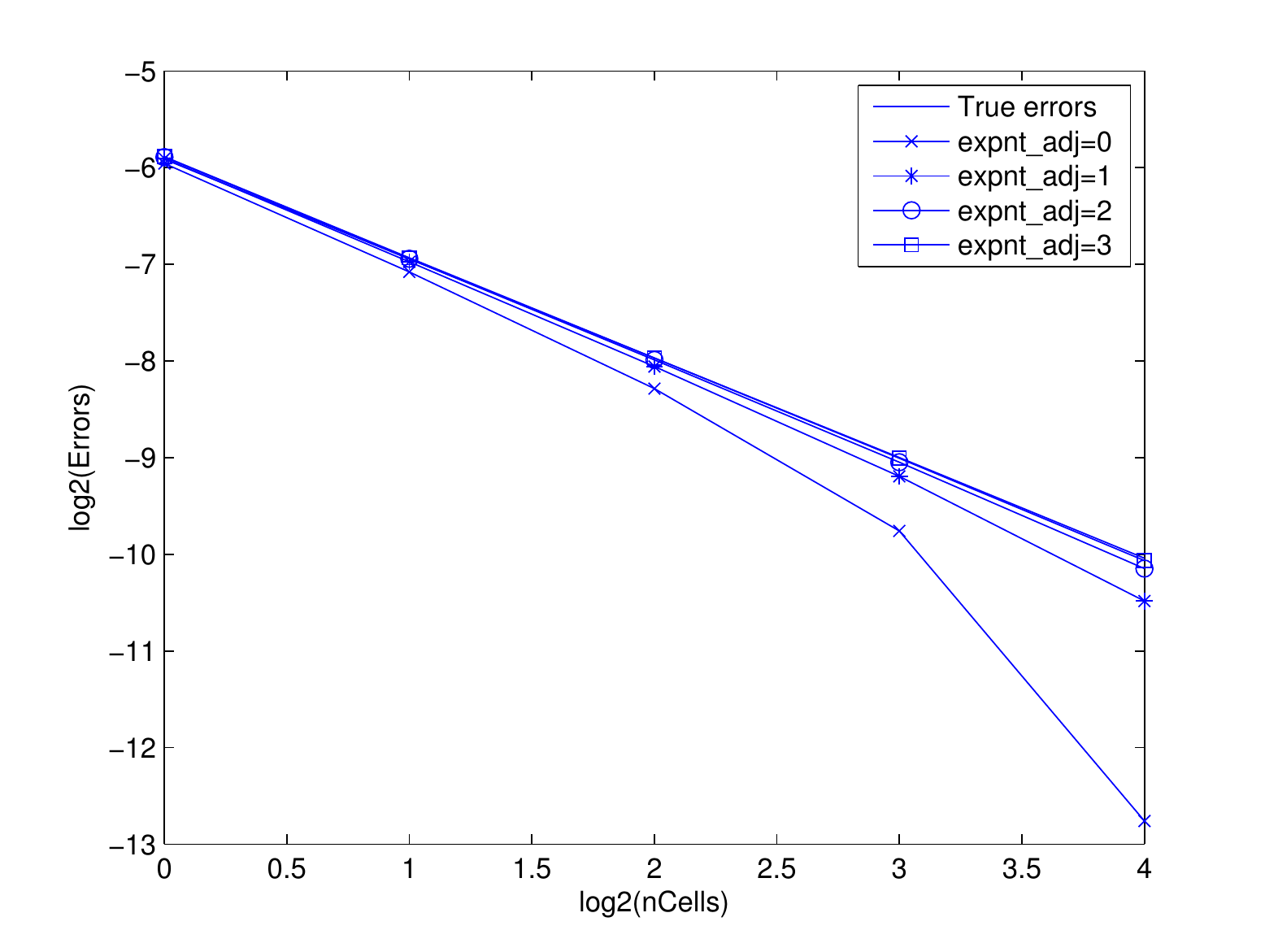} 
   \caption{For the simple kernel
     function $\phi = 1$, the errors in the quantity of interest for different grids for the primal and adjoint approximations. The same first-order method is used for discretizing the primal and adjoint equations.}
   \label{fig:simple1}
   \end{center}
\end{figure}

For the kernel function $\phi$ in \eqref{eq:8}, we first
consider the simple case 
\begin{equation}
\label{eq:10}
   \phi =1. 
\end{equation}
This choice of kernel function emphasizes  the accuracy of the
solution everywhere in the spatial-temporal region $[0,\,1]\times [0,\,T]$.
With this choice of kernel
function, we compute the solution on an array of uniform grids, 
from a coarse one with $M = 20*2^0$ cells to a fine one with
$M = 20*2^4$ cells. Because the solution of \eqref{e1a.1}--\eqref{e1a.3}
is just the sine wave given by \eqref{eq:9}, we 
can compute the {\it true error} in the quantity of interest
$Q(u) = \int_\Omega u\ud x\ud t$. We also compute the solution of the
adjoint equation, using the same first-order upwind method, on  
another array of uniform grids, from a coarse one with
$M_{adj} = 20*2^0$ to a fine one with $M_{adj} =
20*2^3$. With {\it 
  each} of these approximate adjoint solutions, we compute 
an a posteriori estimate of the error in the quantity of 
interest using  \eqref{e1.17}. Then we plot the errors in the quantity
of interest against  
the number of grid cells $M$ in Fig.~\ref{fig:simple1}.
For this simple case, solving the adjoint equation with the same
first-order upwind method on a grid that is at most as fine as the
grid for the primal equation proves adequate. In fact, 
with $M_{adj} = 20*2^3$ (one level coarser than the finest grid
for 
solving the primal equation), the a posteriori error estimates 
are indistinguishable from the true errors. 

However, with a more challenging kernel function, solving 
the adjoint equation with a numerical method of the same order 
on a grid with similar resolution may not be adequate. In such
cases, increasing the grid resolution for the adjoint equation 
can solve this difficulty. A more effective solution strategy is to use a 
higher-order method for the adjoint equation. We demonstrate
these observations with the choice 
\begin{equation}
   \phi = \dfrac{1}{\pi\epsilon^2}\exp(-\dfrac{(x-L/2)^2+(t-T/2)^2}{\epsilon^2}).
   \label{e1a.5}
\end{equation}
This kernel function emphasizes the accuracy of the solution in a
small region of radius $\epsilon$ surrounding the point $(x=L/2,
t=T/2)$ in the space-time computational domain.

\begin{figure}[h!]
   \begin{center}
    {\includegraphics[width=4.5in]{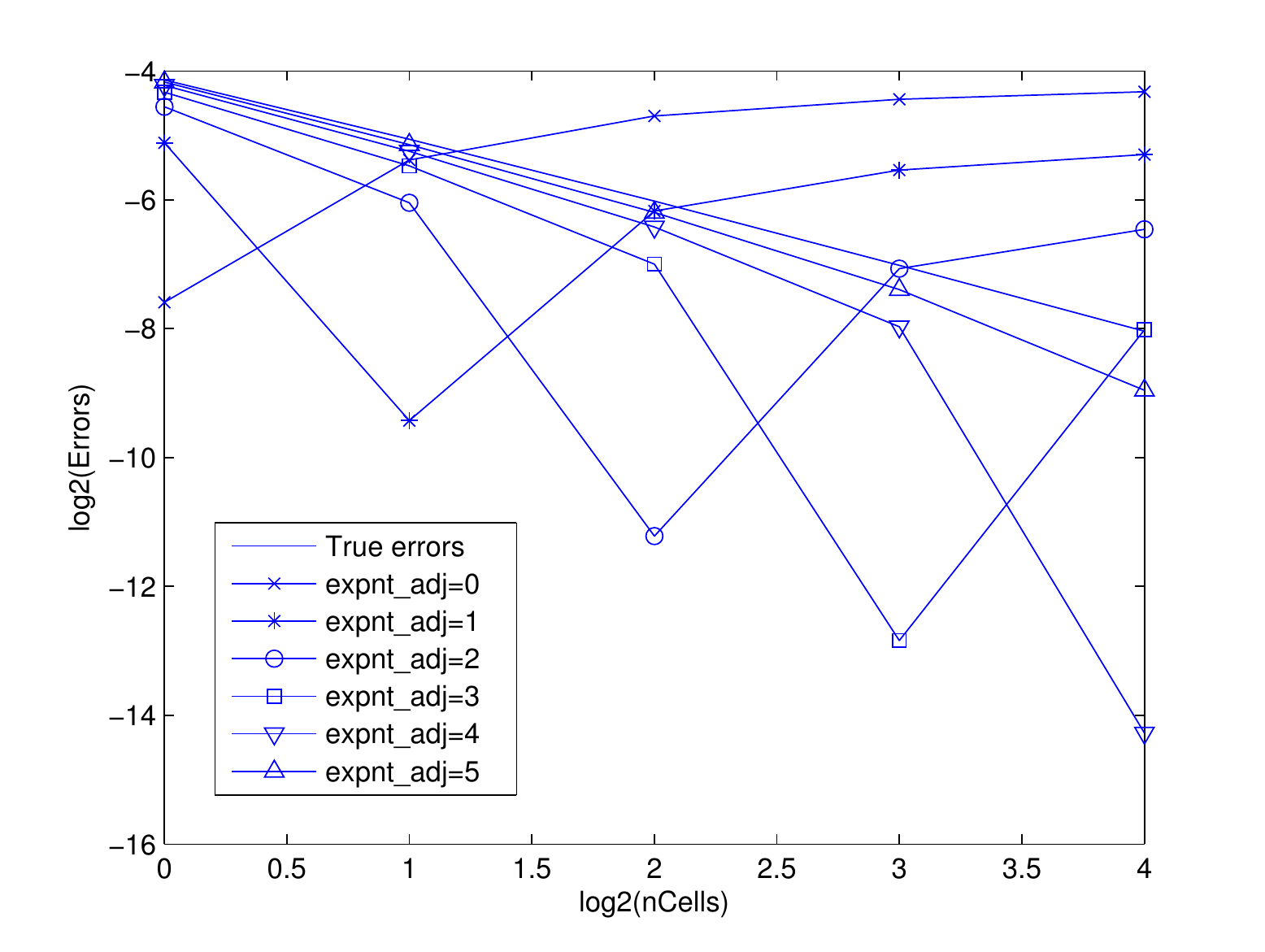}}
   \caption{For the kernel function \eqref{e1a.5}, the errors in the
     quantity of interest for different grids for the primal and
     adjoint approximations. The same first-order method is used for
     discretizing the primal and adjoint equations.} 
   \label{fig:hard1}
      \end{center}
\end{figure}

We first solve the adjoint equation with the same first-order 
upwind method as used for the primal problem on an array of uniform
grids with increasing resolutions. 
The results are presented in Fig.~\ref{fig:hard1}.
What we see is that on coarse grids 
(in this case $M_{adj} = 20*2^q$ with $q=0,\,1,\,2,\,3$), 
the adjoint solutions are not accurate enough to produce
 sensible estimates of the errors in the quantity of 
 interest. We also see that with increasing resolutions for the
 adjoint approximation, 
 the a posteriori estimates of the errors are converging 
 to the true errors. In fact, at $M_{adj} = 20*2^5$, 
 the errors estimates can be regarded as decent approximations
 of the true errors. 

\begin{figure}[h!]
   \begin{center}
   {\includegraphics[width=4.5in]{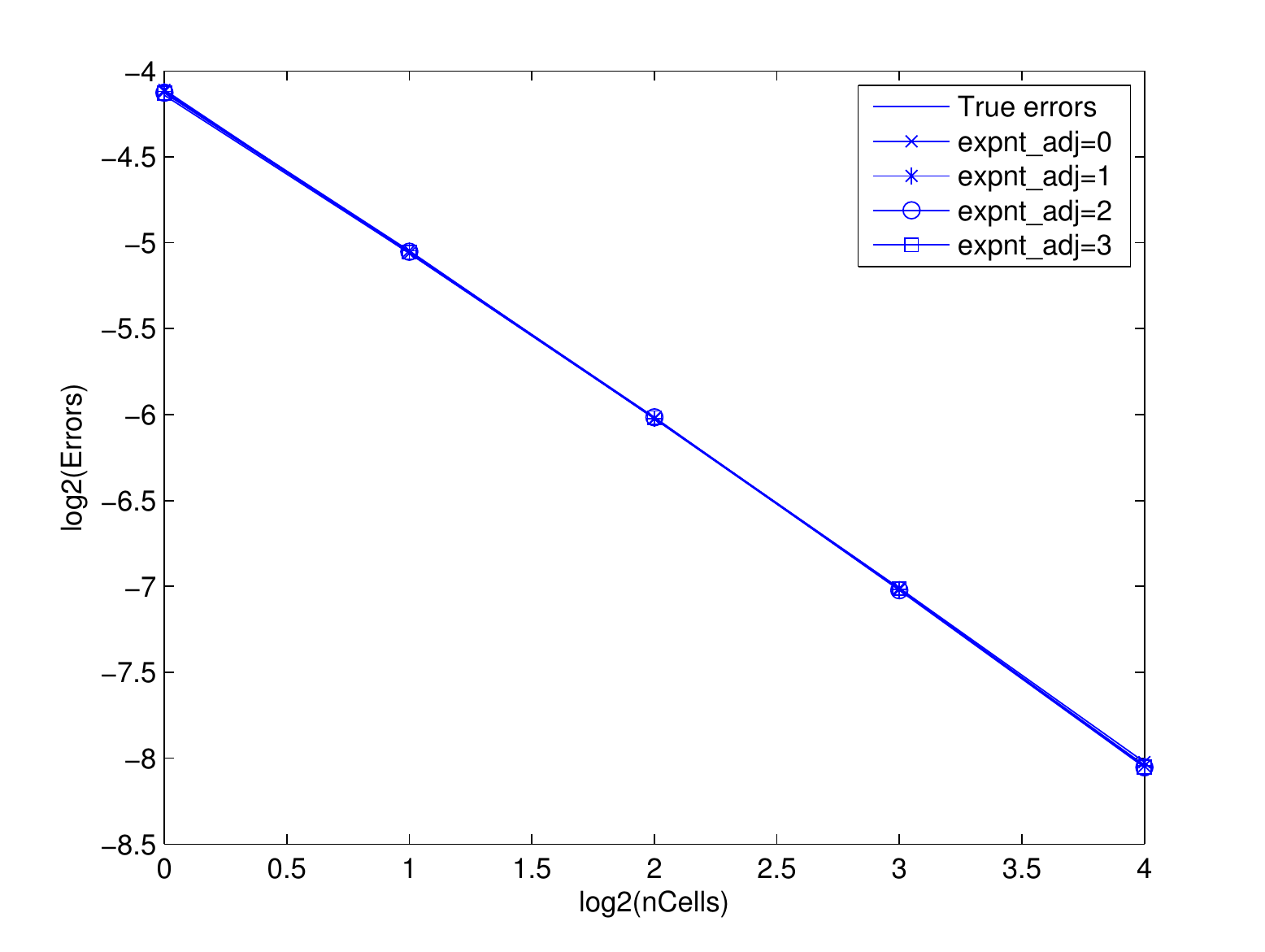}}
   \caption{For the kernel function \eqref{e1a.5}, the errors in the
     quantity of interest for different grids for the primal and
     adjoint approximations. A higher-order method is used for
     discretizing the adjoint equations compared to that used for the
     primal equations.} 
   \label{fig:hard_leap}
      \end{center}
\end{figure}

Next, we change our strategy, and solve the adjoint equation
using the leap-frog method which is second-order accurate 
both in time and space. The results are shown in Fig.~\ref{fig:hard_leap}.
The leap-frog method provides a very efficient solution.
In fact, even on the coarsest grid with just $M_{adj} = 20$, 
the adjoint solution is accurate enough to produce 
error estimates that are indistinguishable from the 
true errors. 

In the Galerkin--type  variational approach towards a posteriori error
estimation for finite volume schemes (\cite{BarLar02, EPPTW09}), it is
required that the solution of the adjoint problem be sought in a
function space larger than the solution space for the primal problem,
due to the orthogonality between the residual and the test function
space. This requirement does not apply to our approach, because
neither of the finite volume schemes for the primal or adjoint
problems is based on variational formulations. In fact, for the cases
\eqref{eq:10} and \eqref{e1a.5} considered in this section, both the
primal and adjoint problems are solved in piecewise constant function
spaces, so to speak. The more complex case \eqref{e1a.5} poses some
challenges because the errors in the adjoint solution starts to
deteriorate the a posteriori estimates of the errors in the quantity
of interest corresponding to this kernel function. Solving the adjoint
problem in a broader function space, e.g.~by implementing a
Godunov-type higher-order finite volume scheme (\cite{BarLar02,
  Lev02}), is certainly a possible strategy. This strategy is not
explored in this work. Instead, we demonstrate through numerical
results that the challenges can also be met, to various extents,
solely by refining the mesh or implementing a plain high order scheme
(e.g.~from a first-order upwind scheme to a second-order leap-frog
scheme). The function spaces for the solutions of the adjoint problems
remain piecewise constant. 

\section{Application to adaptive mesh refinement for the 1D shallow
  water equations}\label{sec:appl-adapt-mesh} 

In this section, we demonstrate the application of the a priori error
estimate, derived in Section \ref{sec:framework}, to adaptive mesh
refinement for finite volume methods. We 
take the system of one dimensional linearized shallow water equations
as an example. We also discuss the challenge and the handling of a
nonlinear quantity of interest. 

\subsection{Types of adaptive mesh refinement strategies}\label{sec:adapt}
\begin{figure}[ht]
   \begin{center}
      \includegraphics[width=4in]{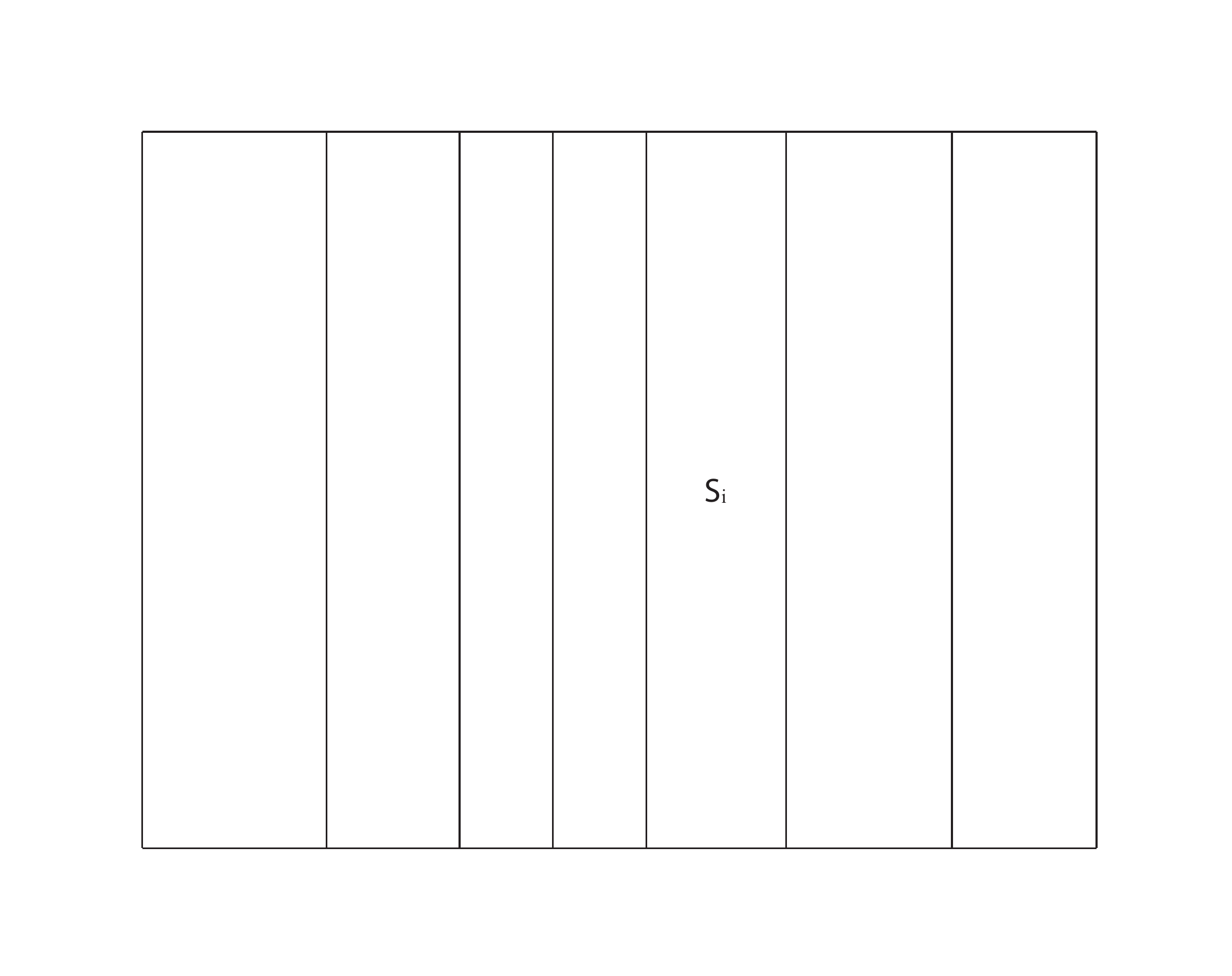}
   \caption{Sketch of a Type 1 adaptively refined mesh.}
\label{fig:type1}
   \end{center}
\end{figure}

\begin{figure}[ht]
   \begin{center}
     {\includegraphics[width=4in]{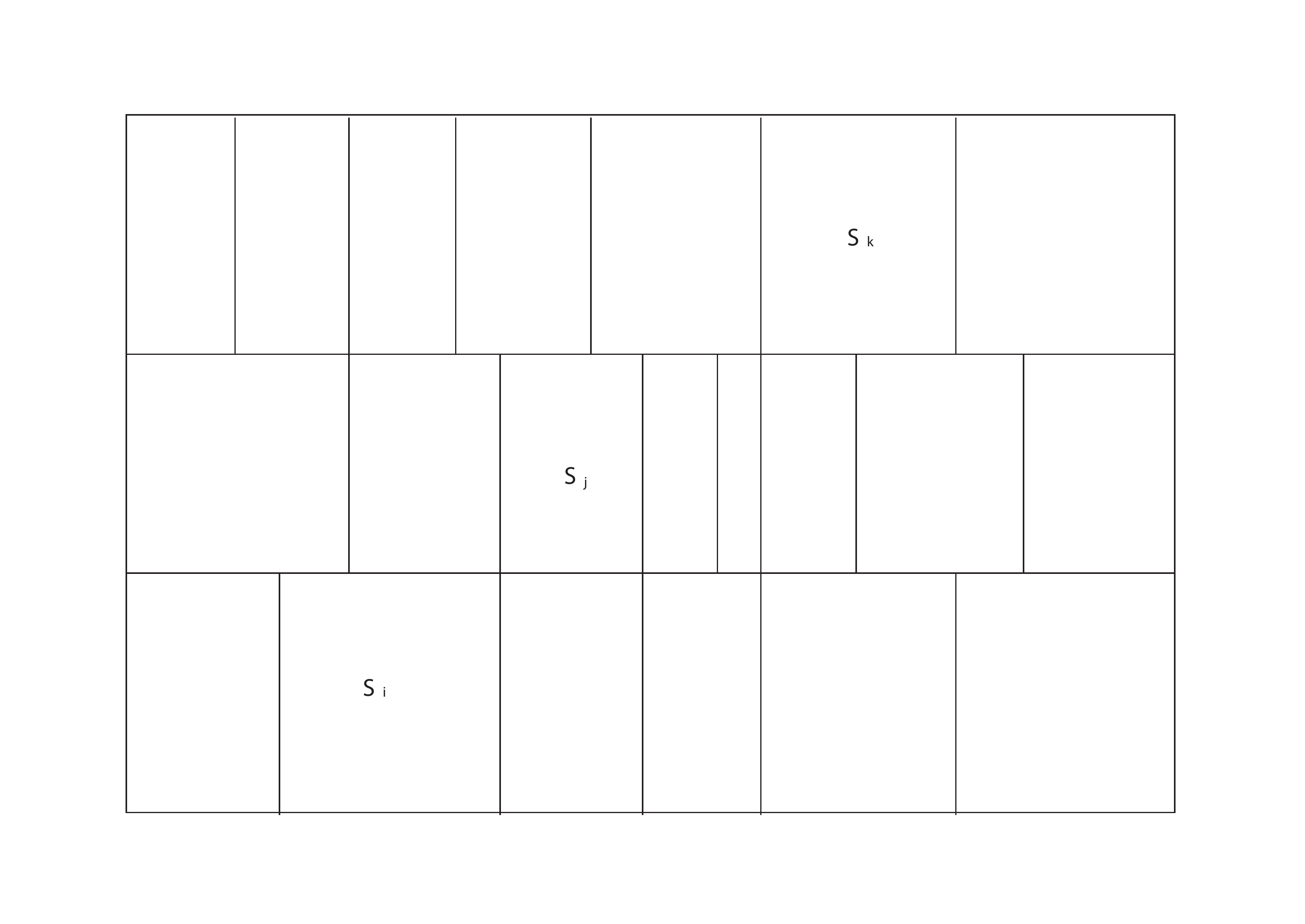}}
   \caption{Sketch of a Type 2 adaptively refined mesh.}
\label{fig:type2}
   \end{center}
\end{figure}

We can identify the following types of adaptive mesh refinement strategies:
\begin{description}
   \item[Type 1.] {\it Time unvarying adaptive mesh refinement.}\\
      The mesh for the spatial dimensions is adaptively refined,
      but remains fixed for the whole simulation period. The
      time steps are non-adaptive, though they could be 
      nonuniform. A generic depiction of the resulting grid is given in
      Fig.~\ref{fig:type1}. 

   \item[Type 2.] {\it Non-incremental time varying adaptive mesh refinement.}\\
      The spatial mesh is adaptively refined for each temporal sub-interval of the
      simulation. The time 
      step is non-adaptive, though
      it could be nonuniform. The mesh refinement is done after each full
      simulation, i.e.,~not incrementally in time. A generic depiction
      of the resulting grid is given in Fig.~\ref{fig:type2}

   \item[Type 3.] {\it Incremental time varying adaptive mesh refinement.}\\
      The spatial mesh is adaptively refined at each sub-interval of the
      simulation period, as the 
      simulation proceeds. The time step is non-adaptive, though
      it could be nonuniform. 

   \item[Type 4.] {\it Incremental time varying adaptive mesh with adaptive
      time steps.}
\end{description}
In this article, we discuss the Type 1 and 2 strategies and
leave the more sophisticated Type 3 and 4 strategies to future work. 

The elegant formula \eqref{eq:6} cannot be applied
directly for adaptive mesh refinements because it does not provide
information about local errors. In order to obtain such information,
we need to break down the whole space-time domain
into slabs, and find a means for calculating the 
error contribution from each slab. We note that  each slab can contain
one or more of the time steps used to for discretization of the
partial differential equations.  
Let $\{t_j\}_{j=0}^N$ be a partition of the whole simulation period
$[0,\,T]$, and for each $j$, let $\{K_{ij}\}_{i=1}^{M_j}$ be a
partition of the spatial domain $(0,\,1)$ during the time period
$[t_{j-1}, t_j]$. Then each slab $S_{ij}$ is given by 
\begin{displaymath}
  S_{ij} = K_{ij}\times t_j.
\end{displaymath}
We note that this partition of the space-time domain can accommodate
the Type 1 adaptive mesh refinement strategy with $N=1$; see
Fig.~\ref{fig:type1}. With $N > 1$, it can accommodate the adaptive
mesh refinement strategies of Type 2 as well as of Type 3 and 4; see
Fig.~\ref{fig:type2}. 

By \eqref{eq:6}, 
\begin{align*}
   Q(u) - Q(\ub^\#) &= (f-L\ub^\#, \tilde\ub)\\
   &= \int_I(f-L\ub^\#)\tilde\ub\ud x\ud t\\
   &= \sum\int_{S_{ij}}(f-L\ub^\#)\tilde\ub\ud x\ud t\\
   &=\sum_{i,j} E_{ij},
\end{align*}
where 
\begin{equation*}
   E_{ij} = \int_{S_{ij}}(f-L\ub^\#)\tilde\ub\ud x\ud t
\end{equation*}
denoting the error contribution from each slab $S_{ij}$.

Below is a simple adaptive mesh refinement algorithm based on
the error contribution from each slab: 
\begin{lstlisting}
   Step 0. Specify the total error tolerance TOL.
   Step 1. Calculate E_ij for each i and j.
   Step 2. Calculate Total_error = sum of E_ij
   Step 3. If Total_error >= TOL, then
              For 1 <= j <= N, 1 <= i <= M_j
                 If |E_ij| >= TOL/M then
                    Refine the grid over the cell K_ij
		 End
              End
	      Goto Step 1.
	   Else if Total_error < TOL, then
	      Exit
	   End
\end{lstlisting}
We note that this adaptive mesh refinement algorithm usually leads to
over refinement because it does not account for the error
cancellations among grid cells. Nevertheless, it is adequate for our
purpose to demonstrate the application of the a posteriori error
estimate \eqref{eq:6} to adaptive mesh refinement for finite volume
methods. For more discussions on adaptive mesh refinement algorithms, see, e.g.~\cite{EEHJ95}.

\subsection{1D shallow water equations}
\label{sec:1d-shallow-water}
We now apply the a posteriori error estimate
\eqref{eq:6} to  adaptive mesh refinement for 
the one-dimensional, linearized shallow water system
\begin{equation}\left\{\begin{aligned} 
       &\dfrac{\p h}{\p t} +  h_x + u_x = 0, \\
       &\dfrac{\p u}{\p t} + 2h_x + u_x = 0.
     \end{aligned}\right.
     \label{e3.1}
\end{equation}
We remark that, among many things, the shallow water system can be
used to model tsunami waves. In \eqref{e3.1}, the system has been
non-dimensionalized, coefficients involving physical quantities have
be replaced by somewhat artificial constants, and the mathematically
non-essential 
Coriolis forcing terms 
have been omitted, all for the sake of simplicity. Despite these heavy
simplifications, the system \eqref{e3.1} still retains some interesting
physical features, e.g. gravity waves (the underlying
mechanism of tsunamis \cite{P87, Maj03}), and therefore is adequate for the purpose of this
section. 

The coefficient matrix of the system \eqref{e3.1}, given by
\begin{equation*}\left(\begin{aligned}
&1 &1\\
&2 &1
\end{aligned}\right),
\end{equation*}
has two eigenvalues
\begin{align*}
  \lambda_+ = 1 + \sqrt{2} > 0\qquad\mbox{and}\qquad
  \lambda_- = 1 - \sqrt{2} < 0.
\end{align*}
The system is sometimes referred to as the {\it subcritical} mode of the
primitive equations due to the the opposite signs of the
eigenvalues. For a discussion on a related problem, see
\cite{CLRTT07}.
The eigenvectors of the coefficient matrix form a
transformation matrix $P$ such that 
$$
P^{-1}\left(\begin{aligned}
&1 & &1 \\
&2 & &1\end{aligned}\right) P = 
\left(\begin{aligned}
&\lambda_+ & &0 \\
&0 & &\lambda_-\end{aligned}\right).
$$
We let 
\begin{equation}
\left(\begin{aligned}
&\xi \\ &\eta\end{aligned}\right) =
P^{-1}\left(\begin{aligned} &h\\ &u\end{aligned}\right).
\label{e3.3}
\end{equation}
Then, \eqref{e3.1} is transformed to 
\begin{equation}
\dfrac{\p }{\p t}\left(\begin{aligned}\xi\\ \eta\end{aligned}\right) +
\left(\begin{aligned} &\lambda_+ &0\\ &0
    &\lambda_-\end{aligned}\right)
\dfrac{\p}{\p x}\left(\begin{aligned} \xi\\ \eta \end{aligned}\right) =
  0. \label{e3.4}
\end{equation}
We prescribe the upwind boundary conditions for $\xi$ and $\eta$ given by
\begin{equation}\left\{\begin{aligned}
& \xi = 0 \qquad \mbox{at $x = 0$},\\
& \eta = 0 \qquad \mbox{at $x = 1$}.
\end{aligned}\right.
\label{e3.5}
\end{equation}
In summary, the primal problem consists of \eqref{e3.1} (or \eqref{e3.4})
and the boundary conditions \eqref{e3.5}.
We let $\Omega = (0,1)\times (0, T)$ and  $\ub$ denote the
vector of unknowns $(h,\,u)$ and let 
\begin{equation*}
L\ub = \left(\begin{aligned}
&h_t + h_x + u_x \\
&u_t + 2h_x + u_x\end{aligned}\right).
\end{equation*}
The domain of the operator $L$ is then defined as 
\begin{equation*}
D(L) = \left\{ \ub \in L^2(\Omega),\,\,\, L\ub \in L^2(\Omega),\,\,\, 
\textrm{and $\ub$ satisfies \eqref{e3.5}}\right\},
\end{equation*}
where $\ub$ is related to $\xi$ and $\eta$ through \eqref{e3.3}.
Therefore, $L$ is an unbounded operator in $L^2(\Omega)$ with domain
$D(L)$.

We shall next define the domain $D(L^*)$ of the adjoint operator. A function $\tilde u$ belongs to $D(L^*)$ if and
only if $\ub \Longrightarrow (L\ub,\,\tilde\ub)$ is a linear
continuous functional on $D(L)$ for the norm of $L^2(\Omega)$. Let
$\ub\in D(L)$ and let $\tilde\ub\in C^\infty(\Omega)$. Then,
\begin{align*}
(L\ub,\,\tilde\ub) &= \int_\Omega (h_t + u_0h_x + h_0u_x)\tilde h +
(u_t + g_0h_x + u_0u_x)\tilde u\\
&\int_0^1h\tilde h\Big|^t_0 - \int_0^Th\tilde h_t +
u\tilde u\Big|_0^T - \int_0^Tu\tilde u_t +
\int_\Omega A\ub_x\cdot \tilde\ub\\
&= \int_0^1(h\tilde h + u\tilde u)\Big|_{t=T} - 
\int_\Omega h\tilde h_t + \int_\Omega AP\left(\begin{aligned}\xi_x\\
\eta_x\end{aligned}\right)\cdot\tilde \ub\\
&= \int_0^1(h\tilde h + u\tilde u)\Big|_{t=T} - 
\int_\Omega h\tilde h_t +\int_\Omega PP^{-1}AP\left(\begin{aligned}\xi_x\\
\eta_x\end{aligned}\right)\cdot\tilde \ub\\
&= \int_0^1(h\tilde h + u\tilde u)\Big|_{t=T} - 
\int_\Omega h\tilde h_t +\int_\Omega \left(\begin{aligned}
&\lambda_+ &0\\
&0 &\lambda_-\end{aligned}\right)\left(\begin{aligned}\xi_x\\
\eta_x\end{aligned}\right)\cdot P^T\tilde \ub\\
&= \int_0^1(h\tilde h + u\tilde u)\Big|_{t=T} - 
\int_\Omega h\tilde h_t +\int_\Omega \left(\begin{aligned}
&\lambda_+ &0\\
&0 &\lambda_-\end{aligned}\right)\left(\begin{aligned}\xi_x\\
\eta_x\end{aligned}\right)\cdot \left(\begin{aligned} \tilde\xi \\
\tilde\eta\end{aligned}\right),
\end{align*}
where 
\begin{equation}
\left(\begin{aligned}\tilde\xi\\ \tilde\eta\end{aligned}\right)
= P^T\tilde\ub.
\label{e3.6}
\end{equation}
Integrating by parts on the space interval we obtain
\begin{align*}
(L\ub,\,\tilde\ub) = &\int_0^1(h\tilde h + u\tilde u)\Big|_{t=T} - 
\int_\Omega h\tilde h_t + \int_0^T\left(\lambda_+\xi\tilde\xi |_{x=1}
- \lambda_-\eta\tilde\eta|_{x=0}\right) - { }\\
\phantom{=} &\int_\Omega\lambda_+\xi\tilde\xi_x + \lambda_-\eta\tilde\eta_x\\
= &\int_0^1(h\tilde h + u\tilde u)\Big|_{t=T} \ud x+
 \int_0^T\left(\lambda_+\xi\tilde\xi |_{x=1}
- \lambda_-\eta\tilde\eta|_{x=0}\right) \ud t - { }\\
\phantom{=} &\int_\Omega h\tilde h_t + u\tilde u_t + \lambda_+\xi\tilde\xi_x
+\lambda_-\eta\tilde\eta_x\\
= &\int_0^1(h\tilde h + u\tilde u)\Big|_{t=T} \ud x+
\int_0^T\left(\lambda_+\xi\tilde\xi |_{x=1}
- \lambda_-\eta\tilde\eta|_{x=0}\right) \ud t - { }\\
\phantom{=} &\int_\Omega h\tilde h_t + u\tilde u_t + \left(\begin{aligned}\xi\\
\eta\end{aligned}\right) \left(\begin{aligned}
&\lambda_+ &0\\
&0 &\lambda_-\end{aligned}\right) \left(\begin{aligned}\tilde\xi_x\\
\tilde\eta_x\end{aligned}\right)\\
= &\int_0^1(h\tilde h + u\tilde u)\Big|_{t=T} \ud x+
\int_0^T\left(\lambda_+\xi\tilde\xi |_{x=1}
- \lambda_-\eta\tilde\eta|_{x=0}\right) \ud t - { }\\
\phantom{=} &\int_\Omega h\tilde h_t + u\tilde u_t + 
\ub\cdot A^T\dfrac{\p}{\p x}\tilde\ub\\
= &\int_0^1(h\tilde h + u\tilde u)\Big|_{t=T} \ud x+
 \int_0^T\left(\lambda_+\xi\tilde\xi |_{x=1}
- \lambda_-\eta\tilde\eta|_{x=0}\right) \ud t - { }\\
\phantom{=} &\int_\Omega \ub\cdot\left(\dfrac{\p\tilde\ub}{\p t} +
A^T\dfrac{\p}{\p x}\tilde\ub\right).
\end{align*}
For $\ub \longrightarrow (L\ub,\,\tilde\ub)$ to be continuous for the
$L^2$ norm, it is necessary that 
\begin{equation}\left\{\begin{aligned}
&\tilde\ub = 0\qquad\textrm{ at } t = T,\\
&\tilde\xi = 0\qquad\textrm{ at } x = 1,\\
&\tilde\eta = 0\qquad\textrm{ at } x = 0,\\
&\dfrac{\p\tilde\ub}{\p t} + A^T\dfrac{\p}{\p x}\tilde\ub\in
L^2(\Omega).
\end{aligned}\right.
\label{e3.7}
\end{equation}
In the above, $(\tilde\xi,\,\tilde\eta)$ is the transformation of
$\tilde\ub$ defined in \eqref{e3.6}.
Therefore, we define the domain for the adjoint operator as
\begin{equation}
D(L^*) = \left\{ \tilde\ub\in L^2(\Omega) \,\,\,|\,\,\, \tilde\ub \textrm{
    satisfies } \eqref{e3.7}\right\}.
\label{e3.8}
\end{equation}
For each $\tilde\ub\in D(L^*)$, 
$$
L^* \tilde\ub = -\dfrac{\p\tilde\ub}{\p t} - A^T\dfrac{\p}{\p
  x}\tilde\ub.
$$
The adjoint problem can be stated as follows:\\
\indent{\em for every $\bs{\phi}\in L^2(\Omega)$, find $\tilde\ub\in
D(L^*)$ such that }
$$
L^*\tilde\ub = \bs{\phi}.
$$
The well-posedness of the adjoint problem can be established by the
semigroup theory.

    We consider the {\em nonlinear} quantity of interest:
    \begin{equation}
    \label{eq:12}
	Q(h,\,u) = \dfrac{1}{2}\int_\Omega hu^2
    \end{equation}
    which is the time integral of the kinetic energy.
    Let $h^\#$ and $u^\#$ be the numerical approximations to $h$ and $u$, respectively.
    Then, the error in the quantity of interest can be computed as
    \begin{align*}
	Q(h,\,u) - Q(h^\#,\,u^\#) &= \int_\Omega \dfrac{1}{2}hu^2 - \int_\Omega \dfrac{1}{2}h^\#{u^\#}^2 \\
	&=\int_\Omega(h-h^\#)\dfrac{1}{2}u^2 + \int_\Omega (u-u^\#)\dfrac{1}{2}h^\#(u+u^\#).
    \end{align*}
Therefore, the kernel function is given by
\begin{equation}
\label{eq:1}
   \phi = \left(\dfrac{1}{2}u^2, \, \dfrac{1}{2}h^\#(u+u^\#)\right).
\end{equation}
We note that the kernel function $\phi$ involves the unknown function
$u$, an issue inevitable for nonlinear quantities of interest. In
practice, the true values of the unknown functions are not available,
and therefore have to be approximated. To the best of our knowledge,
there is yet no rigorous theory guding the choice of the
approximations. An obvious option, which is also what is usually taken
in the literature, is to replace the unknowns by their numerical
approximations. Because the goal here is to evaluate the applicability
of our a posteriori error estimation approach to adaptive mesh
refinement, we content ourselves with this option, and leave the more
fundamental questions to future endeavor. Thus in what follows, we
replace $u$ by $u^\#$ in \eqref{eq:1}.

\begin{figure}[h!]
  \centering
  \includegraphics[width=4in]{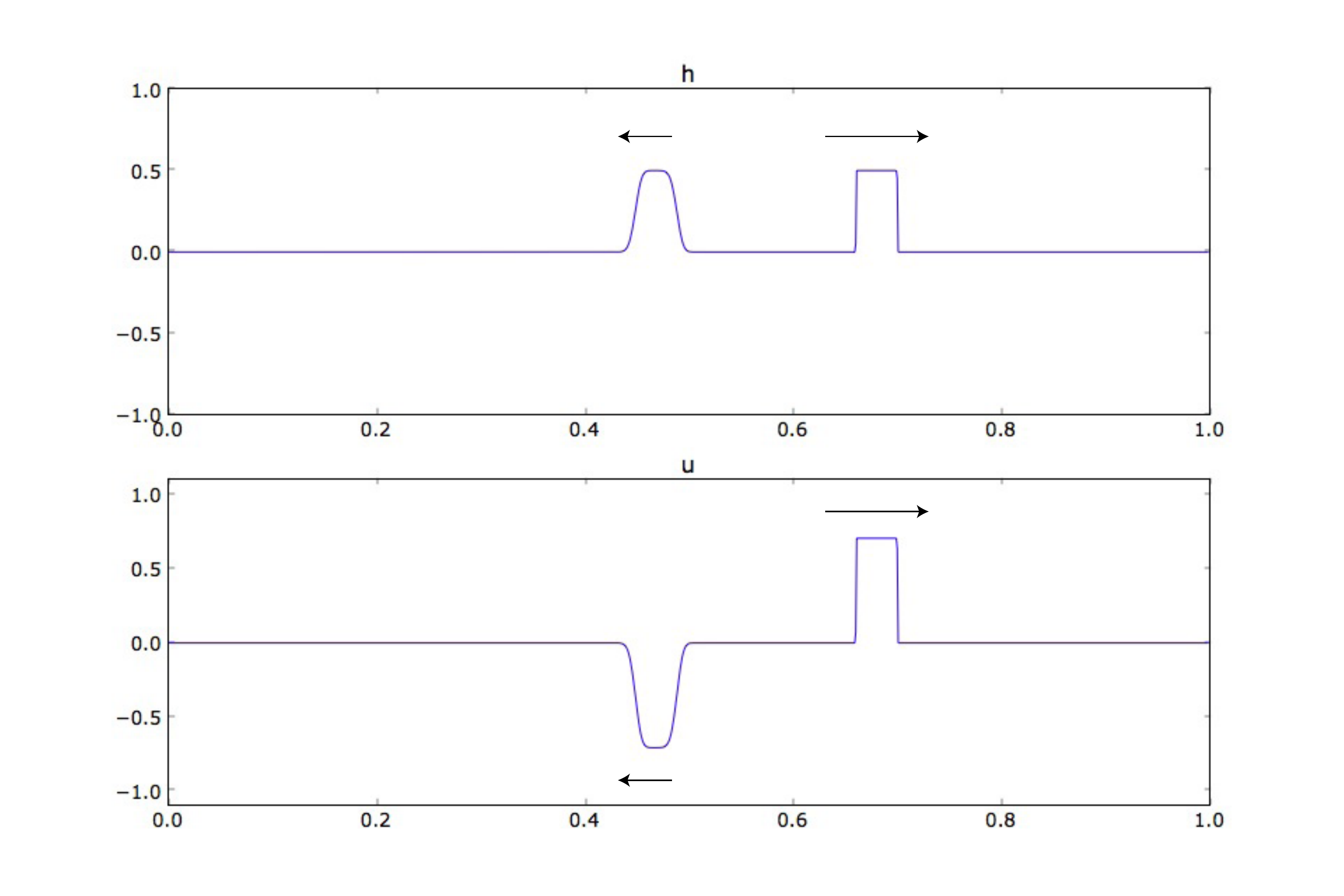}
  \caption{Snapshot of the solution of \eqref{e3.1}. The arrows point
    to the wave travel directions; longer arrows
    indicate greater speeds, but not proportionally.}
  \label{fig:snapshot}
\end{figure}

For the system \eqref{e3.1},  we specify as the initial conditions
    \begin{equation*}
	h(x,0)= h_0(x) = \left\{\begin{aligned}
	    & 1\qquad |x-\frac{1}{2}| < \epsilon,\\
	    & 0 \qquad \textrm{ elsewhere},
	\end{aligned}\right.
    \end{equation*}
and
    \begin{equation*}
	u(x,0) = 0.
    \end{equation*}
The set of initial conditions represents a flow at rest with a bulk of
fluid artificially raised above the surface, reminiscent of the sudden
flow elevation caused by an earthquake under the sea. When the flow is allowed
to evolve freely, the bulk of fluid will split into two smaller wave
packets of equal size and travel in opposite directions with
different speeds. See Fig.~\ref{fig:snapshot} for a snapshot of the
solutions $h$ and $u$.
The problem is numerically challenging due to
the discontinuous data and solutions and to 
the wave packets moving in different directions at different 
speeds.
    
\begin{figure}[h!]
   \begin{center}
\includegraphics[width=4in]{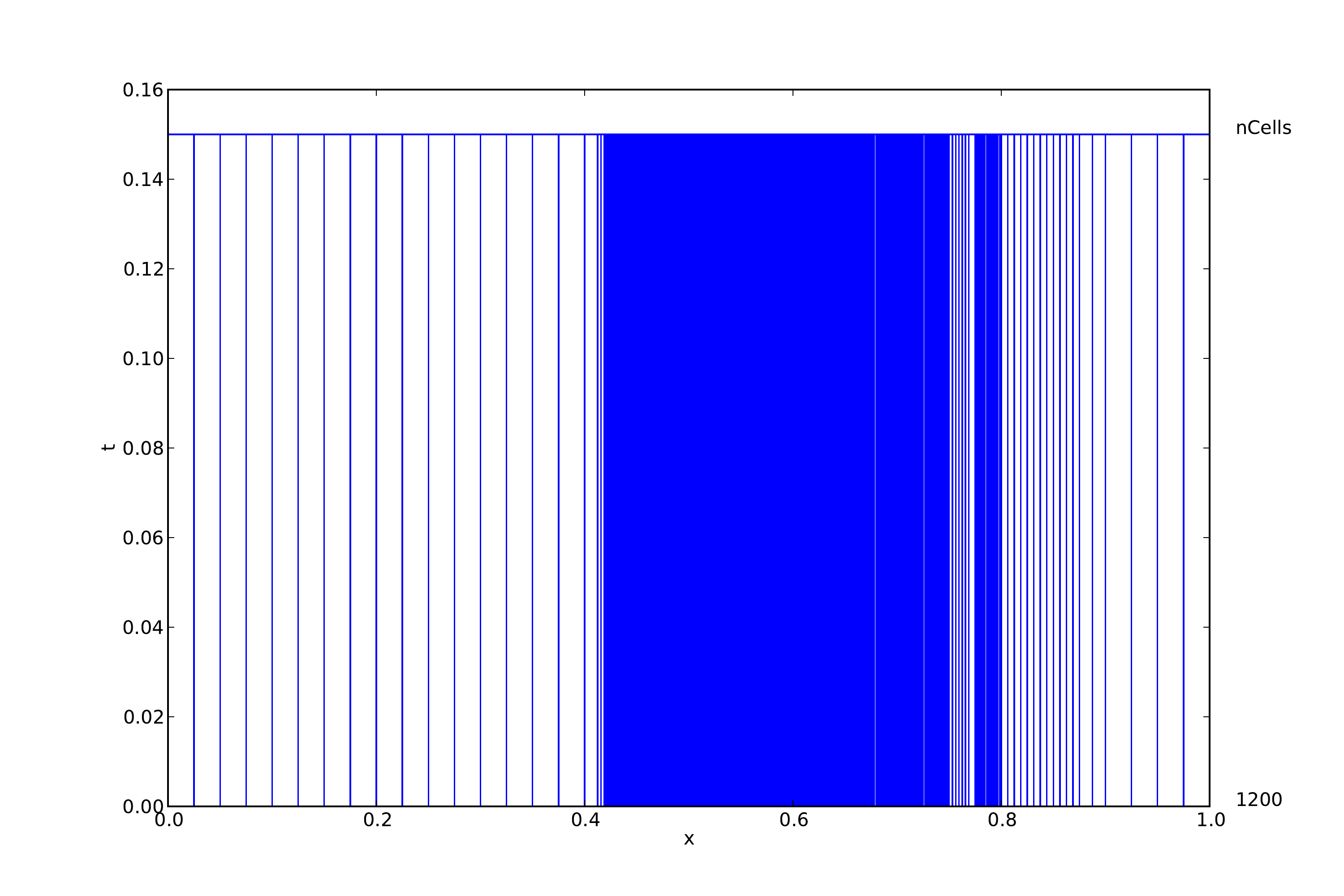}
  \caption{The final grid generated by the Type 1 adaptive mesh
    refinement strategy.}
   \label{fig:grid-type1}
      \end{center}
\end{figure}

\begin{figure}[ht]
   \begin{center}
  \includegraphics[width=4in]{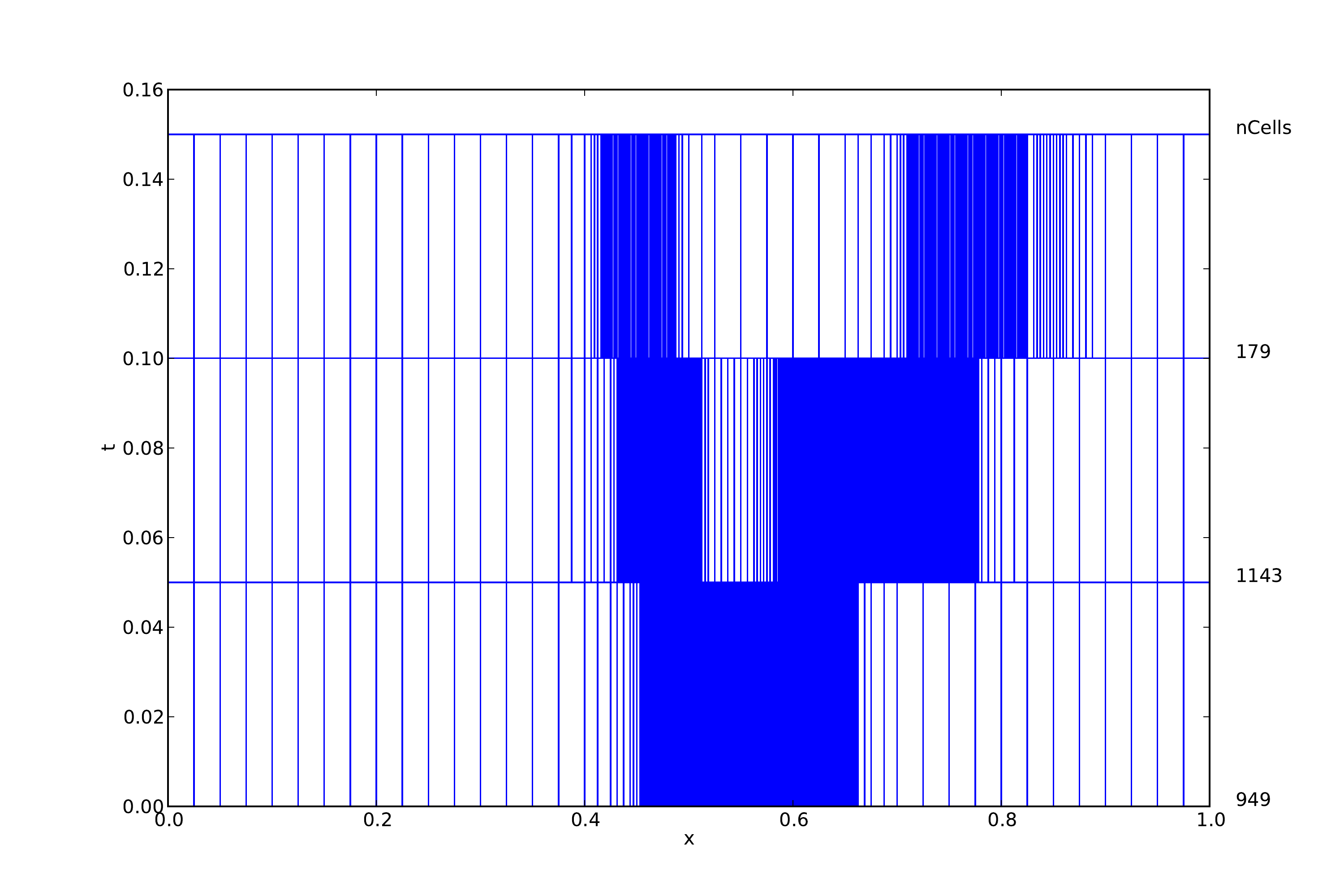}
   \caption{The final grid generated by the Type 2 adaptive mesh
     refinement strategy}
  \label{fig:grid-type2}
     \end{center}
\end{figure}

We experiment with Type 1 and Type 2 mesh refinement
strategies and compare the results with that of the uniform mesh
refinement strategy. All strategies are applied with the same
tolerance goal 
\begin{equation*}
TOL = 4.0\times 10^{-4}
\end{equation*}
for the error. The uniform mesh refinement strategy requires
2,560 cells to reach the tolerance goal. The Type 1 strategy,
i.e., time-unvarying adaptive mesh refinement, requires 1,200 cells
to reach the same tolerance goal. The final grid is plotted in
Fig.~\ref{fig:grid-type1}. We see that the grid is intensely
refined over the travel ranges of both of the wave packets. The savings
in number of cells come from the ``quiet'' regions where there are no
wave activities. The refined region is shifted towards the right end
because the rightward wave is traveling at a greater speed and thus
has a longer travel range.  We then apply the Type 2 strategy with three temporal
sub-intervals over the whole simulation period $[0,T]$. The resulting grid that meets the tolerance
is given in Fig.~\ref{fig:grid-type2} and has 
949 cells for the first sub-interval, 1143 for the second, and 179 for
the last, with an average of 757 cells for the whole simulation
period. The saving in the number of cells stems from the fact that
when the wave packets are far apart, the surrounding regions for them
can be refined separately and thus the unnecessary refinement in the
middle is avoided. 

\section{Concluding remarks}\label{sec:conclusion}
In this article, we present a framework for goal-oriented a posteriori
error estimation for finite volume methods. The formulation of the a
posteriori error estimate is independent of the exact form  of the
methods and therefore can be applied to arbitrary
finite volume methods. In this framework, it is not required that the
adjoint equations be solved in a function space larger than that for
the forward equation, due to the fact that finite volume methods, generally
speaking, are not built on Galerkin orthogonality. 

To demonstrate the validity of the a posteriori error estimate, we
conduct numerical experiments with the one-dimensional linear
transport equation. The primal equation is solved using the
first-order, upwind finite volume method. The overall conclusion from
these experiments is that the more accurate the adjoint solution is,
the more accurate the a posteriori error estimate is. Roughly
speaking, there are two scenarios one has to deal with in practice. If
the quantity of interest is defined by a simple kernel function, then
solving the adjoint equation with a numerical method of the same order or accuracy
and with a
mesh with resolution similar as that for the primal
equation may prove adequate. If the quantity of interest is more
complex, then the adjoint equation may have to be solved on a finer
mesh or to be solved by a high-order accurate method. We have found that the
latter approach is often more efficient. 

An application of the a posteriori error estimate to adaptive mesh
refinement is also presented. The one-dimensional linearized
shallow-water equations are taken as an example. The test case
involves two wave packets traveling in opposite directions with
different speeds. This case is numerically very challenging. The a
posteriori error estimate is found to be effective at guiding various
adaptive mesh refinement strategies that lead to grids that are
dynamically refined according to wave activities. 

We can identify several directions or future work. The impact of the lack of accuracy in the adjoint solution on
the accuracy of the a posteriori error estimate is only experimentally
explored in this work. More rigorous analysis is warranted for this
issue. A related issue is the impact of the substitution of the
primal solution $\ub$ by its approximation $\ub^\#$ in the kernel
function $\phi$. We did not touch upon this issue here, but it is very
important and inevitable in cases involving nonlinear quantities of
interest.   

Another direction future research is to apply the framework of a
posteriori error estimation, laid out in this article, to regional
climate modeling, which is the primary motivation for the current
work. An emerging approach towards climate modeling is to use one
global grid over the whole sphere, with local refinements over regions
of interest, and with smooth transitions between coarse and fine
regions (\cite{RJGJDS}). In our opinion, the current mesh refinement
strategy used in this approach is quite rudimentary in that it only
refines over regions that are directly of interest. Goal-oriented error
estimation is clearly needed to guide a more sensible mesh refinement
strategy, and thus to control the computational errors with regard to
the quantities of interest.  

\section*{Acknowledgment}
The author (QC) thank Du Pham for helpful comments and Varis Carey for 
helpful discussions. This work is supported by  the
Department of Energy grant number DE-SC0002624. 

\bibliographystyle{siam}
\bibliography{amr}

\end{document}